\documentclass[12pt, a4paper]{article}
\usepackage{amsmath,amssymb,amsthm}
\usepackage[final]{pdfpages}

\def\veps{\varepsilon}

\newcommand{\beq}{\begin{equation}}
\newcommand{\eeq}{\end{equation}}
\newcommand{\beqa}{\begin{eqnarray}}
\newcommand{\eeqa}{\end{eqnarray}}

\title{ A sharpening of Tusn\'ady's inequality}

\author{
Jen\H{o} Reiczigel$^1$
\and
L\'{\i}dia Rejt\H{o}$^{2,3}$
\and
G\'abor Tusn\'ady$^3$
}

\begin{document}
\maketitle
\footnotetext[1]{  Szent Istv\'an University, Department of Biomathematics and Informatics, Faculty of Veterinary Science, Budapest,  Hungary}
\footnotetext[2]{University of Delaware, Statistics
Program, FREC, CANR, Newark, Delaware, USA}
\footnotetext[3]{Alfr\'ed  R\'enyi
Mathematical Institute of the Hungarian Academy of Sciences,
Budapest, Hungary}

\begin{abstract}
Let  ~$\veps_1, \ldots, \veps_m$ be i.i.d. random variables with
 $$P(\veps_i=1)= P(\veps_i= -1)=1/2,$$
and $X_m = \sum_{i=1}^m \veps_i.$  Let $Y_m $ be a normal random variable with the same first two
moments as that of $X_m.$  There is a uniquely determined function $\Psi_m$ such that the distribution
of  $\Psi_m(Y_m)$ equals to the distribution of $X_m$. Tusn\'ady's inequality states that
$$ \mid \Psi_m(Y_m) - Y_m \mid  \leq \frac{Y_m^2}{m}+1.$$ Here we propose a sharpened version of this
inequality.

\end{abstract}

\section{Conjecture}
Let ~$\veps_1, \ldots, \veps_m$ be i.i.d. random variables with
$$P(\veps_i=1)= P(\veps_i= -1)=1/2,$$ and $X_m = \sum_{i=1}^m \veps_i.$
Let $Y_m $ be a normal random variable with the same first two moments
as that of $X_m.$ Using quantile transformation we can see that there is
a uniquely determined function $\Psi_m$ such that the distribution of
$\Psi_m(Y_m)$ equals to the distribution of $X_m$.  The central limit
theorem implies that the function $\Psi_m$ is close to the identity for
large $m$.  A sharp inequality of Tusn\'ady \cite{TG} raised certain
interest in the literature
(\cite{BM},\cite{Cas},\cite{Pol},\cite{CsM},\cite{CsR},\cite{Dud},\cite{Dud2},\cite{Maj},\cite{MasD},\cite{Mas},\cite{Pol2},\cite{Zh}).

Let us define the function $f$ on the interval $(0,1)$ as
$$
f(x)= \sqrt{(1+x)\log(1+x)+(1-x)\log(1-x)},
$$
set $f(0)=0, f(1)=\sqrt{\log(4)}.$
Let us put
$$
x_{k,m} = {\frac{k-{\frac{m}{2}}}{\frac{m}{2}}}
$$
for positive even integers $m$ with $k$ such that
$m/2 < k \leq m,$ and set
$$
p_{k,m} = P(X_m \geq 2 k-m) = 2^{-m} \sum^m_{i=k}{m \choose i}.
$$
Let us define the function $Q$ on the reals as
$$
Q(x) = {\frac {1}{\sqrt{2\pi}} \int^\infty_x e^{-u^2/2}du}.
$$
With those ingredients our conjecture states that
$$
Q({\sqrt{m}}f(x_{k,m}))<p_{k,m}<Q({\sqrt{m}}f(x_{k-1,m}))
$$
holds true for ${\frac{m}{2}}<k\leq m$. Or more sharply
\beq \label{eq:CON}
2(k-1)-{\frac{m}{2}}+0.8964 < m f^{-1} (Q^{-1}(p_{k,m})/\sqrt m) <
2(k-1)-{\frac{m}{2}}+1.0000
\eeq
holds true with pessimal parameters $m=k=10.$
It implies that Tusn\'ady's inequality is sharpened to
$$
\left | \Psi_m(Y_m) - mf^{-1}\left({\frac{Y_m}{m}}\right) \right | < 1.1036.
$$

\section{Generalization}

For an arbitrary random variable $X$
let us consider the  function on reals
$$
 R(t) = E e^{tX}
$$
restricting ourselves for distributions having finite
momentum generators. Next we define
$$
\psi(t) = {\frac {R'(t)}{R(t)}},
$$
$$
\alpha(x) = t \quad {\text{iff}} \quad \psi(t) = x,
$$
$$
\rho(x) = R(\alpha(x)) \exp(-x \alpha(x)).
$$
The probability $P(\sum^m_{i=1} X_i \geq m x)$ is approximately
$\rho(x)^{-m}$ if $x>E X$. The function $\rho$ depends on the distribution
of $X$, it is the Chernoff function of $X$.
Let us denote the Chernoff function  of the distribution
$F$ of $X$ by $\rho_F$, and the corresponding function for
standard normal by $\rho_G.$ The quantile transformation between
the partial sums of distribution $F$ with Gaussian ones resemble us
to the equation
$$
\rho_F (x) = \rho_G (y)
$$
having the property that it gives sharp values for any $m.$
Perhaps the error term is bounded with a bound depending on
the distribution of $X.$ For the case symmetrical binomial
distribution the error term might be as small as that
the quantile curve jumps over its limiting function:
it is the informal explanation of our  conjecture.

\section{Numerical Illustration}
The function $\Psi_m$ is shown in Figure 1. called ``step" for $m=50$
with a rescaling for random variables $$\xi_m = \frac{X_m}{m}, \quad
\eta_m= \frac{Y_m}{m}.$$ The function $f$ is called ``limit", for the
sequence of step functions goes to $f$ after rescaling.  The conjecture
comes from the observation that the limit function crosses all steps
near to their middle.  Let us introduce the blow up error term
$$
\Delta_{k,m} = 10 \left( 2 k-1-m \; f^{-1}  \left(   \frac{1}{\sqrt{m}} \; Q^{-1} \left( \sum^m_{i=k} {m \choose i} 2^{-m} \right) \right) \right),
$$
for $0 < k \leq m/2.$ In Figure 1. it is labelled as "Delta''.  With
these notations (\ref{eq:CON}) is equivalent with $ 0< \Delta_{k,m} <
1.036 .$ These error terms are shown in Figure 2. for $2 \leq m \leq
1000.$ Figure 2. prompts the conjecture that even these curves are
convergent.  We are a bit perplexed:  even the inequality
$0<\Delta_{1,2}<1.036$ means that $Q(0.723359)<0.25<Q(0.6435214).$ How
can we prove such an inequality theoretically?

\bigskip\noindent
{\sc Acknowledgement:} We thank to Peter Harremo\"es for pointing out a mistake in
the earlier version of the paper.

\newpage
\section{Appendix}
\noindent
{\bf{ {\tt{R}}- program of Figures 1 and 2.}

\begin{verbatim}
Q=function(p) -qnorm(p)
G=function(x) ((1+x)*log(1+x)+(1-x)*log(1-x))**0.5
Ginv=function(u) {
GG=function(x) G(x)-u
uniroot(GG,c(0,1),f.lower=-u,f.upper=log(4)^.5-u,tol=10^-100)
}

m=50; k=m/2
sum=0; divisor=2**m; bin=
xx=c(1:k+1); yy=c(1:k+1); zz=c(1:k+1);

for (i in 1:k-1){
sum=sum+bin
x=(m-2*i)/m
y=Q(sum/divisor)/(m**.5)
b=Ginv(y)$root
yy[i+1]=y; xx[i+1]=x
bin=(m-i)*bin/(i+1)
zz[i+1]=10*(m-2*i-1-m*b)}
xx[k+1]=0; yy[k+1]=0; zz[k+1]=0
kerx=c(0,1.25); kery=c(0,1.15)
plot(kerx, kery, type="n",xlab="eta", ylab="xi",
main="Figure1. Quantile transform, its limit and blownup error, m=50")

for (i in 1:k){
bb=seq(from=yy[i+1], to=yy[i], by=0.01)
cc=bb*0+1; cc=cc*xx[i+1]
points(bb,cc,type="l", col="blue", lwd=2)}
cc=seq(from=0, to=0.999, by=0.001)
bb=((1+cc)*log(1+cc)+(1-cc)*log(1-cc))**0.5
points(bb,cc, type="l", col="red", lwd=2)
points(yy,zz, type="l", col="green", lwd=2)
legend(locator(1),c("Limit","Step","Delta"),
lty=c(1,1,1),
col=c("red","blue","green"))


kerx=c(0,1.25); kery=c(0,1.15)
plot(kerx, kery, type="n",xlab="eta", ylab="Delta",
     main="Figure 2. The blownup error")

for (k in 1:500){m=2*k;
sum=0; divisor=2**m; bin=1
yy=c(1:k+1); zz=c(1:k+1);
for (i in 1:k-1){
sum=sum+bin
y=Q(sum/divisor)/(m**.5)
b=Ginv(y)$root
yy[i+1]=y;
bin=(m-i)*bin/(i+1)
zz[i+1]=10*(m-2*i-1-m*b)}
yy[k+1]=0; zz[k+1]=0
if (k<100) clr="red" else
if (k<200) clr="blue" else
if (k<300) clr="purple" else
if (k<400) clr="gray" else clr="green"
points(yy,zz, type="l", col=clr)}
legend(locator(1),c("0<m <= 200","200<m<=400","400<m<=600",
"600<m<=800","800<m<=1000"),
lty=c(1,1,1,1,1),
col=c("red","blue","purple","gray","green"))
\end{verbatim}

\newpage
\begin{figure*}[t]
\centering
\includegraphics{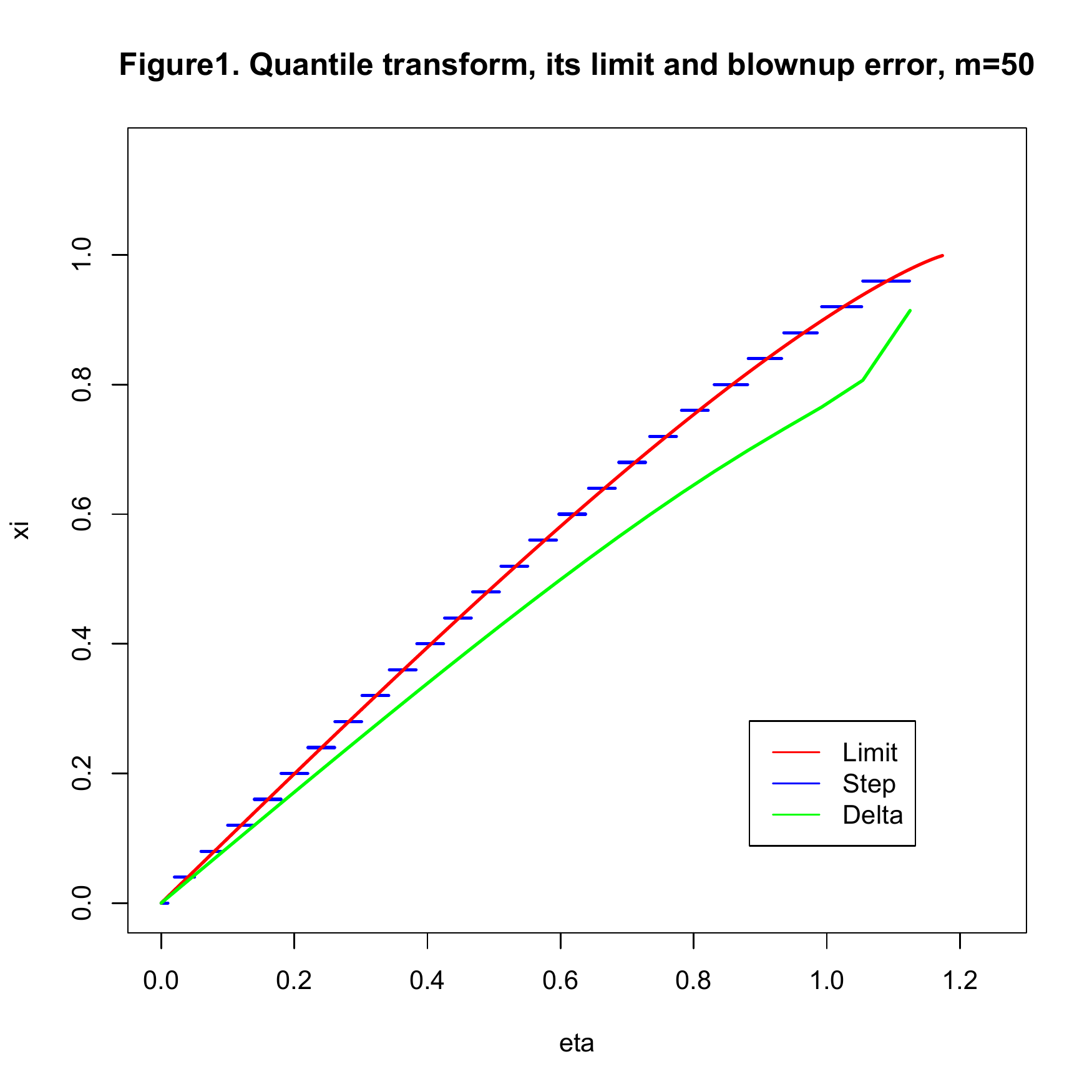}
\end{figure*}

\newpage
\begin{figure*}[t]
\centering
\includegraphics{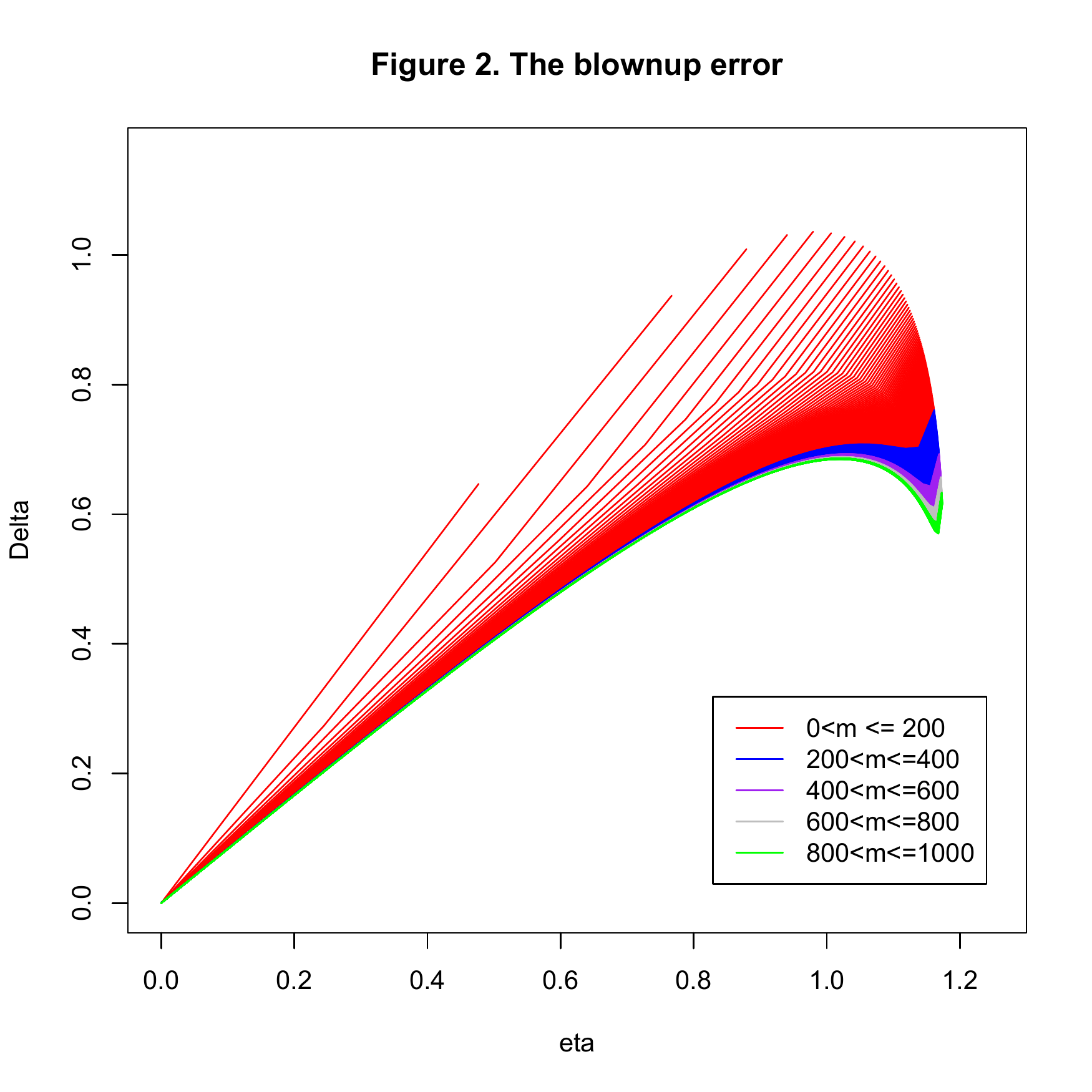}
\end{figure*}


\bigskip\bigskip
\begin{thebibliography}{99}

\bibitem{BM} Bretagnolle, J. and Massart, P. (1989).  Hungarian
constructions from the nonasymptotic viewpoint, {\it The Annals of
Probability} {\bf 17}, 239--256.

\bibitem{Cas} Castelle, N. (2009). Improvement of two Hungarian
bivariate theorems, Manuscript

\bibitem{Pol} Carter, A. and Pollard, D. (2004).  Tusn\'ady's
inequality revisited, {\it The Annals of Statistics} {\bf 32},
2731--2741.

\bibitem{CsM} Cs\"org\H o, M. (2007). A glimpse of the KMT (1975):
approximation of empirical processes by Brownian bridges via quantiles,
Acta Sci. Math. (Szeged) {\bf 73}, 349--366

\bibitem{CsR} Cs\"org\H o, M. and R\'ev\'esz, P. (1981).  {\it Strong
Approximations in Probability and Statistics},  Academic Press, New
York.

\bibitem{Dud} Dudley, R. M. (1999). {\it {Notes on empirical processes,}}
Lecture notes for a course given at Aarhus Univ., August 24, 1999.

\bibitem{Dud2} Dudley, R. M. (2008). On the quantile transformation for
asymmetric binomial and hypergeometric distributions, Manuscript

\bibitem{Maj} Major, P. (1999). {\it {The approximation of the empirical
distribution function,}} Technical report, Alfr\'ed R\'enyi Mathematical Institute of the
Hungarian Academy of Sciences.  Notes available from
{\mbox{http://www.renyi.hu/ \~ .major/probability/empir.html}}, August 3, 1999.

\bibitem{MasD} Mason, D. M. (2001).  Notes on the KMT Brownian bridge
approximation to the uniform empirical process,  In {\it Asymptotic
Methods in Probability and Statistics with Applications}, (N.
Balakrishnan, I. A. Ibragimov and V. B. Nevzorov, eds.) 351--369.
Birkh\"auser, Boston.

\bibitem{Mas} Massart, P. (2002).  Tusn\'ady's lemma, 24 years later,
{\it Ann.  Inst.  H. Poincar\'e Probab.  Statist.} {\bf 38}, 991--1007.

\bibitem{Pol2} Pollard, D. (2002). {\it A user's guide to
measure theoretic probability,} Cambridge University Press

\bibitem{TG} Tusn\'ady, G. (1977).  {\it Study of statistical
hypotheses}, Dissertation for Habilitation, Hungarian Academy of
Sciences, Budapest.

\bibitem{Zh} Zhou, H. H. (2006). A note on quantile coupling inequalities
and their applications, Manuscript

\end{thebibliography}
\end{document}